# Design of Hybrid Regrouping PSO-GA based Sub-optimal Networked Control System with Random Packet Losses

Indranil Pan[a], Saptarshi Das[a,b,]*

a) Department of Power Engineering, Jadavpur University, Salt Lake Campus, LB-8, Sector 3, Kolkata-700098, India.
b) Communications, Signal Processing and Control Group, School of Electronics and Computer Science, University of Southampton, Southampton SO17 1BJ, United Kingdom.

Authors' emails:
indranil.jj@student.iitd.ac.in, indranil@pe.jusl.ac.in (I. Pan)
saptarshi@pe.jusl.ac.in, s.das@soton.ac.uk (S. Das*)

**Abstract:**

In this paper, a new approach has been presented to design sub-optimal state feedback regulators over Networked Control Systems (NCS) with random packet losses. The optimal regulator gains, producing guaranteed stability are designed with the nominal discrete time model of a plant using Lyapunov technique which produces a few set of Bilinear Matrix Inequalities (BMIs). In order to reduce the computational complexity of the BMIs, a Genetic Algorithm (GA) based approach coupled with the standard interior point methods for LMIs has been adopted. A Regrouping Particle Swarm Optimization (RegPSO) based method is then employed to optimally choose the weighting matrices for the state feedback regulator design that gets passed through the GA based stability checking criteria i.e. the BMIs. This hybrid optimization methodology put forward in this paper not only reduces the computational difficulty of the feasibility checking condition for optimum stabilizing gain selection but also minimizes other time domain performance criteria like expected value of the set-point tracking error with optimum weight selection based LQR design for the nominal system.

# 1. Introduction:

Networked Control Systems deals with the enforcement of control policies over an unreliable real time network and has been the subject of active research in the recent past. This increasing trend of implementing control systems over a communication network is partly attributed to the cheap availability of off-the-shelf communication hardware and the various advantages of using a shared communication medium, viz. reduced wiring, modularity, etc. (Tipsuwan and Chow 2003). Thus NCS is advantageous both from the economic as well as implementation point of view and has been adopted in diverse fields ranging from industry automation to space applications.

However NCS is not a panacea and has a few pressing issues which need to be addressed before actual implementation on a real plant. Due to the unreliable nature of the communication network, the packets may get delayed due to queuing at the buffers and suffer stochastic delays at each sampling instant. In some cases due to buffer overflows or bit errors the packets are lost which is referred to as packet dropout. These packet losses severely hamper the control system performance and thus must be accounted for at the design stage itself. It is not unnatural for a plant to become unstable in the presence of packet dropouts and thus exponential stability of such systems must be guaranteed.

In (Zhang et al. 2001), the control system with the network is modeled as an asynchronous dynamical system (ADS) and stability conditions for the system with packet dropouts is ensured through the formulation of a set of bilinear matrix inequalities. However the BMI problem is non-convex and is computationally difficult to solve. This makes the computation of the stabilizing state feedback gains difficult and the conventional approach is to resort to random search methods as in (Zhang et al. 2001). However in process control applications, other constraints like faster rise time and settling time, zero steady state error etc. need to be incorporated in the design phase itself while simultaneously ensuring stability of the system in the presence of packet losses. Thus a more tractable solution scheme is required which honors the stability constraints of the BMIs and at the same time give acceptable time domain performance in the presence of packet losses. In this paper a hybrid Regrouping PSO-GA based design methodology has been adopted to achieve this objective. In (Longo et al. 2012; Longo et al. 2011) global optimization techniques like GA and PSO have been used for NCS problems as well. However, the stochastic algorithms are used to address the issues of scheduling in NCS. In the present paper, the hybridization of PSO-GA algorithms along with standard LMI solver has been used for ensuring stability of the closed loop control system under arbitrary packet losses. In this paper, the proposed hybrid algorithm also minimizes the expected global minima of the random set-point tracking error which occurs due to the stochastic nature of packet loss in NCS for fixed state feedback controller, while also satisfying analytical Lyapunov stability theorem i.e. the BMI criterion using hybridization of standard LMI solver and GA optimizer.

To the best of authors' knowledge such hybrid algorithm has not been previously applied for solving BMI problems often encountered in NCS stabilization. Along with ensuring guaranteed Lyapunov stability in the presence of random packet losses, we optimized for performance (using time domain performance criterion for state feedback control law) with near optimal deviation of the state variables and control signal due introduction of the discrete time LQR formulation. This new design methodology incorporates a two stage global optimizer for efficient solution of the problem which opens a wide range of possibilities to include other performance criteria in the control design process. This makes the problem to be applied in a more practical setting than a theoretical framework as a real control system design must be designed for a wide variety of performance objectives and not just Lyapunov stability for stabilization problem over communication network.

The rest of the paper is organized as follows. Section 2 discusses the theoretical formulation for handling packet losses in NCS. Section 3 discusses the hybrid PSO-GA based suboptimal state-feedback regulator design methodology and Section 4 gives an illustrative example to show the validity of the proposed method. The paper ends in Section 5 with the conclusions followed by the references.

# 2. Theoretical formulation for handling packet losses in NCS

## 2.1. Formulation of the stability criterion as a BMI problem

The NCS with random packet losses can be considered as an asynchronous dynamical system with the following structure (Zhang et al. 2001) as shown in Fig. 1. ADSs like hybrid systems incorporates continuous modes (like differential or difference equation) and discrete modes (discrete event driven finite automata) as well. The derivation and formulation of the BMI problem as in (Zhang et al. 2001), is discussed next.

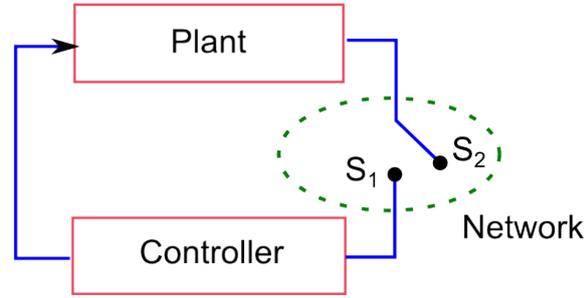

Fig. 1. Schematic diagram showing the plant and the state feedback controller with the network represented by a switch on one side.

Let us take the plant is governed by the following difference equation

$$x(k+1) = Gx(k) + Hu(k) \tag{1}$$

where, $G$ and $H$ are discrete-time system matrices. The system is controlled by the state feedback law over NCS.

$$u(k) = -K\bar{x}(k) \tag{2}$$

Here, $x(k)$ and $\bar{x}(k)$ are the states of the system before and after the switch respectively in Fig. 1. Now if $r_1, r_2, \cdots, r_N$ be the rates representing the fraction of time of occurrence of each discrete state, then $\sum_{i=1}^{N} r_i = 1$. If there exist a Lyapunov function $V(x(k))$ and scalars $a_1, a_2, \cdots, a_N$ corresponding to each rate such that

$$a_1^{r_1}, a_2^{r_2}, \cdots, a_N^{r_N} > a > 1 \tag{3}$$

and $V(x(k+1)) - V(x(k)) \leq (a_s^{-2} - 1)V(x(k)), s = 1, 2, \cdots, N.$ (4)

then the ADS remains exponentially stable with decay rate greater than $a$. Now, if the discrete state dynamics is given by $x((k+1)h) = \tilde{\Phi}_s x(kh), s = 1, 2, \cdots, N$, the search for the quadratic Lyapunov function of type $V(x(kh)) = x^T(kh)\tilde{P}x(kh)$ and the scalars $a_1, a_2, \cdots, a_N$ can be cast into a bilinear matrix inequality or BMI problem. Thus equation (3) and (4) reduces to

$$a_1^r a_2^{1-r} > 1,$$ (5)

$$\tilde{\Phi}_1^T \tilde{P} \tilde{\Phi}_1 \leq a_1^{-2} \tilde{P},$$ (6)

$$\tilde{\Phi}_2^T \tilde{P} \tilde{\Phi}_2 \leq a_2^{-2} \tilde{P}$$ (7)

for $s = 1, 2$.

When the switch is in position $S_1$.

$$\tilde{\Phi}_1 = \begin{bmatrix} G & -HK \\ G & -HK \end{bmatrix};$$ (8)

and when the switch is in position $S_2$,

$$\tilde{\Phi}_2 = \begin{bmatrix} G & -HK \\ 0 & I \end{bmatrix}.$$ (9)

In (5)-(7), for a specified gain matrix $K$, $\tilde{\Phi}_1$ and $\tilde{\Phi}_2$ become constant. Thus, the solution to these equations reduces to the finding of the scalars $a_1, a_2$ and a positive semi-definite matrix $\tilde{P}$.

### 2.2. Computational complexity of the BMI problem

Linear Matrix Inequalities or LMIs are convex in nature and efficient interior point methods (Vandenberghe and Balakrishnan 1997) exist to solve these easily. The set of Equations (6)-(7) are BMI in nature due to the multiplication term of the solution variables $a_1, a_2$ with $\tilde{P}$ which is also another solution variable itself. This makes the system non-convex in nature. Also BMI problems are known to be NP hard (Toker and Ozbay 1995). The class P refers to problems which are solvable in polynomial time. A lot of research is going on in the theory of computational complexity regarding the solution techniques of a NP-hard problem and it is still an open question. It is generally assumed that an NP-hard problem cannot have a polynomial time solution in the worst case. NP-hard is a not a characteristic of any particular algorithm but of the problem itself and practical algorithms for solving these type of problems typically involve approximations or heuristics (VanAntwerp and Braatz 2000). Branch and bound algorithm as in (VanAntwerp and Braatz 2000) is one of the ways of solving BMIs but requires tight upper and lower bounds for the objective function which is difficult to obtain in many cases.

There exist a variety of methods to solve BMI problems, commonly encountered in various control related problems. Although there are no efficient interior point methods for BMI unlike those of the LMI, there are many instances where depending on the specific structure of the formulated control problem a mathematical simplification may be made. The problem can then be solved by seeking the solution to an equivalent set of LMIs. However these techniques are not universal (i.e. they cannot be applied to all BMI problems), and often such analytical mathematical transformations are difficult to obtain and are not trivial. For example, in (Lo and Lin 2004), the authors note that the problem is non-convex and NP hard. Their control problem is such that it can be cast into the form of parameterized linear matrix inequality (PLMI) technique as proposed in (Tuan et al. 2001). Therefore in spite of the problem being a BMI one, it is possible to solve their particular problem of robust control for fuzzy systems using their proposed two-step method using the PLMI to obtain relaxation criterion and then using standard LMI solvers. There are many other instances where such problem specific simplifications may be made and the BMI problem can be reduced to a set of more computationally tractable set of LMIs. In (Kim and Kim 2002) for example, the stability conditions are derived in terms of BMIs. To solve this set of BMIs an iterative LMI approach (ILMI) is proposed and additional mathematical manipulations are done to relax the constraints and solve it using standard solvers. It is however to be noted that this is a kind of local search algorithm and within a multimodal solution space and therefore it might converge to local minima depending on the initial conditions. This aspect is highlighted further in (Kanev et al. 2004), where it is stated that most of the existing local search approaches are computationally fast but depending on the initial conditions, these might not converge to a global solution. The D-K iteration for μ-synthesis (Doyle 1983), alternating semi-definite programs (SDP) method (Fukuda and Kojima 2001), dual iteration method (Iwasaki 1999), are not guaranteed to converge to local solutions as shown in (Goh et al. 1995)(Fukuda and Kojima 2001)(Yamada and Hara 1998). Hence in (Kanev et al. 2004) a local BMI optimization is proposed using an iterative technique for finding out the globally optimal robust output feedback controller. It is to be noted that though this paper overcomes the problem of global optimization, it is limited to a class of BMI problems (robust controllers for linear systems) and hence cannot be directly incorporated for our case. Other global optimization algorithms are variations of the branch and bound technique (Tuan and Apkarian 2000)(VanAntwerp and Braatz 2000).

In our case the control problem is more complicated as the system is not a linear system as described in all the previous literatures. Since the problem at hand is a switched system, there are two ways to approach it in the light of the previously published works. The first is to try to mathematically propose some transformation and try to frame the present BMI problem in such a way so that it can be solved by standard LMI solvers. This would give guaranteed convergence and an upper bound on the runtime convergence of the algorithm. The second approach is to use hybrid optimization involving evolutionary or swarm based techniques coupled with standard LMI solvers as has been done in the present case. This does not give any upper bound on the runtime of the algorithm or guaranteed convergence, unlike that of the first case but this method is generic and is not limited to the present problem at hand. Thus it can be easily extended to other control design problems without necessitating any complicated mathematical manipulations as required in the other cases. Application of hybrid swarm, evolutionary and LMI based techniques can be found in other control applications as well for example in active suspension system (Kong et al. 2012).

A comparison with respect to the previously published methods has not been done here, as mathematical transforms for converting a BMI to an equivalent set of LMIs have not yet been investigated for the case of switched systems with LQR formulation and robust set-point tracking requirement. Formulating this mathematically would be another investigation in its own right and would be a digression from the main theme of the present study. Therefore the comparison is restricted to the same class of hybrid evolutionary and swarm algorithms using different variants, namely PSO-GA and GA-GA in order to efficiently solve the BMI problem for the present NCS problem.

## 3. Discrete optimal control for NCS

### 3.1. GA based BMI solving technique

In the present approach, genetic algorithm is used to solve the BMIs to find out the stability of the system in the presence of random packet losses. GA has been used to solve matrix inequalities as in (Sekaj and Vesely 2005) where each element of the solution matrix is considered as a decision variable. This technique has been shown to outperform the standard deterministic V-K iteration method in (El Ghaoui and Balakrishnan 1994). However as the dimension of the solution matrix increases the number of variables of GA would increase and it is difficult to find solutions for multimodal functions in this multidimensional search space. Hence an alternative method is adopted in the present paper as described below.

It is well known that GA is a stochastic optimization algorithm inspired from Darwin's theory of evolution. The GA starts with an initial population of randomly chosen solution variables (represented as chromosomes). A fitness function evaluates the fitness of each individual in each generation. These individuals undergo reproduction, crossover and mutation to give rise to newer individuals in the next generation. The solution is thus iteratively refined until the desired fitness value is reached or the maximum number of generations is exceeded. The individuals are ranked according to their fitness value and the higher ranked individuals have a greater probability of going into the next generation by the process of reproduction. Crossover refers to information exchange in a probabilistic fashion among two parent individuals to give rise to a child. In mutation a small part of the solution vector is randomly altered and the new vector goes to the next generation. There is also a parameter called elite count which dictates the number of fittest individuals in the present generation that is definitely copied over to the next generation. Generally the elite count is taken as a small fraction of the total population as otherwise the solution converges prematurely due to dominance of the initially found fitter individuals and the exploration of the solution space at later generations is hindered. In our case, the population is comprised of 20 individuals and the elite count is taken as 2. The remaining individuals undergo crossover and mutation depending on the Crossover Ratio (CR) and the Mutation Ratio (MR) respectively. The CR is chosen to be 0.8 and the MR is chosen to be 0.2.

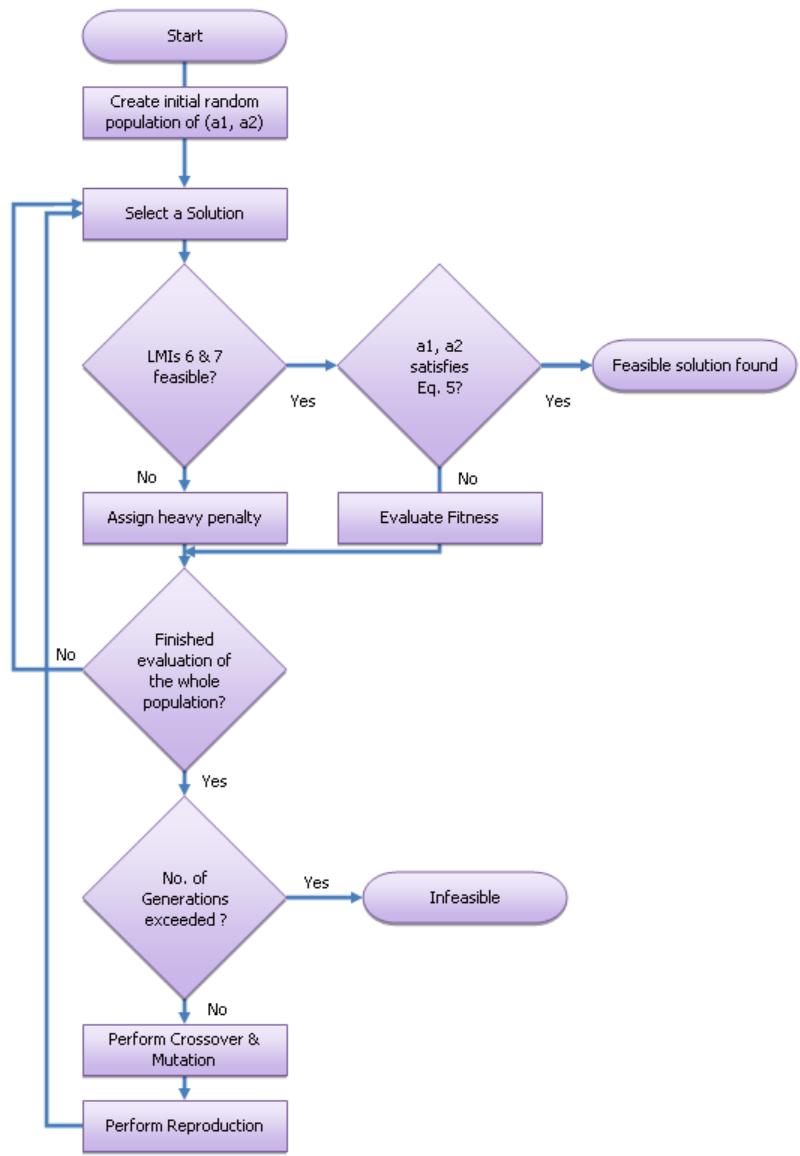

Fig. 2. Flowchart for solution of BMI problem with Genetic Algorithm.

The flowchart for the solution of the GA based BMI problem is shown in Fig. 2. The network is chosen to have 70% packet transmission ($r = 0.7$) similar to (Zhang et al. 2001). The scalars $a_1, a_2$ are chosen to be solution variables of GA. Thus in the fitness function evaluation phase, (6) and (7) are solved using standard interior point method (VanAntwerp and Braatz 2000), since for specified values of $a_1, a_2$ the problem reduces to a LMI equation instead of a BMI one. If a feasible solution is found, (5) is then tested to evaluate feasibility. If not found feasible, the fitness function assigns a value to the solution vector depending on the degree of infeasibility. Equation (5) is modified to form the fitness function as given by the following.

$$f = 1 - a_1^r a_2^{1-r} \tag{10}$$

If the value of this fitness function becomes negative, then a feasible solution is found. Thus the solution is iterated until a feasible solution is found or the maximum number of generations is exceeded. It is to be noted that the algorithm stops on the first occurrence of a feasible solution. Thus all the intermediate solutions are infeasible.

### 3.2. Optimal stabilizing gain selection with Regrouping PSO

In this subsection the basic philosophy of discrete optimal control is first introduced (Ogata 1995). For the discrete system governed by (1), the task is to design an optimal state feedback regulator which minimizes the infinite horizon quadratic optimal cost

$$\tilde{J} = \frac{1}{2} \sum_{k=0}^{\infty} \left[ x^T(k) Q x(k) + u^T(x) R u(x) \right] \tag{11}$$

Minimization of the quadratic cost given in (11) leads to the solution of the Discrete Algebraic Riccati Equation (DARE) given by (12)

$$P = Q + G^T P G - G^T P H \left( R + H^T P H \right)^{-1} H^T P G \tag{12}$$

In (12), $Q, R$ are the weighting matrices and $P$ is the positive definite solution of the Riccati equation (12). Matrix $P$ produces the optimal state-feedback gain matrix $K$ which minimizes the quadratic cost function (11) using the following relation

$$K = \left( R + H^T P H \right)^{-1} H^T P G \tag{13}$$

Thus the optimal control law is given by

$$\begin{aligned} u(k) &= -Kx(k) \\ &= -\left( R + H^T P H \right)^{-1} H^T P G x(k) \end{aligned} \tag{14}$$

But in the presence of the packet losses due to the switch in Fig. 1, optimality can not be preserved as reported in equation (2). Thus a discrete time Linear Quadratic Regulator (LQR) based state feedback controller can neither guarantee the stability nor optimal time domain performance for random packet drop-outs in the NCS. In order to guarantee the stability with a nominally chosen optimal regulator the GA based BMI solution technique is first reckoned. If the state feedback regulator stabilizes the system with a certain probability of packet losses in all combinations, then the weighting matrices of the discrete Riccati equation is chosen optimally with respect to another performance criterion which affects the output of the process as

$$J = \sum_{k=1}^{\infty}\left[k \cdot |e(k)|\right] = \sum_{k=1}^{\infty}\left[k \cdot |r(k) - y(k)|\right] \tag{15}$$

Here, $\{e, r, y\}$ denotes the tracking error, reference input and system response respectively. It is clear that the performance index $J$ in (15) is equivalent to the Integral of Time multiplied Absolute Error (ITAE) criterion in continuous time controller design. The inclusion of the absolute error term reduces the peak overshoot and the steady state error. The multiplication of the time term penalizes the error more at the later stages and hence ensures faster rise and settling times. Also, for a guess value of the state feedback controller gains within the optimization algorithm, due to the stochastic nature of the packet-drop phenomena, the system output $y(t)$ becomes a random variable. Therefore, the objective function (15) signifies expected or average tracking performance for various possible realizations of the packet loss. The calculation of the objective function thus involves multiple-time simulation of the same NCS with the same state feedback controller gains and taking their expected value for minimization.

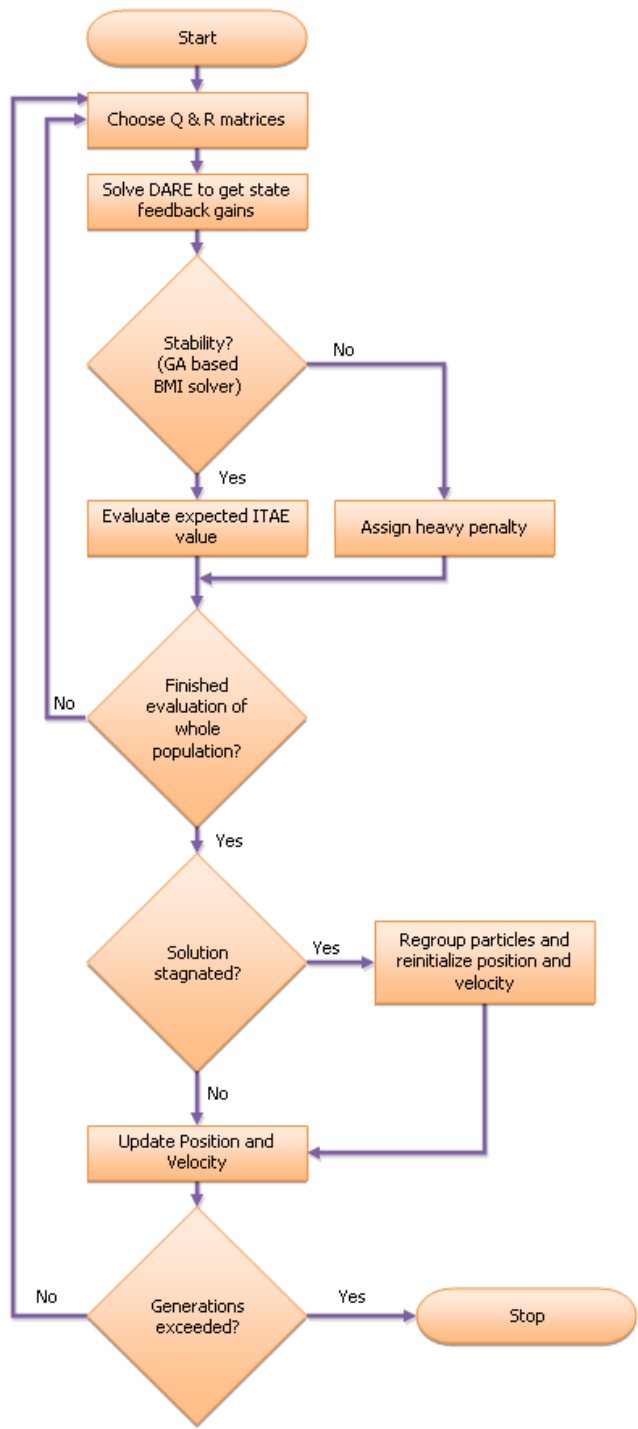

Fig. 3. Flowchart of the Hybrid Regrouping PSO-GA based optimization algorithm for stabilizing gain selection.

The average of the cost function (15) for mutiple runs of the same NCS controller is minimized with Regrouping PSO algorithm while optimally choosing the weighting matrices for the discrete LQR problem for the nominal system without packet losses thus producing sub-optimal controller gains that guarantees the stability using the BMI condition and also gives satisfactory closed loop time response. The overall solution technique for the hybrid regrouping PSO-GA based sub-optimal controller design has been shown in Fig. 3. A heavy penalty is incorporated to discourage search with unstable solutions over the iterations for similar optimization based controller design like in (Zamani et al. 2009). The output from the interior point based LMI solver is just a yes or no for stability, i.e. the LMI solver incorporated in the inner GA loop can only say whether it is stable or not (a binary value) and does not give a continuous value for the degree of stability (like that provided by magnitude of the unstable eigen-values for LTI systems). Hence it is not possible to incorporate a penalty function due to the nature of the problem, though it is well known that the penalty function approach might have been a more expedient method. The Regrouping PSO has certain advantages over the deterministic optimization techniques and even the genetic algorithm which are discussed next along with its working philosophy.

The classical *gbest* PSO is a swarm based stochastic optimization algorithm which does not need any gradient based information for the minimization of an objective function (Kennedy and Eberhart 1995). In *gbest* PSO the particles are initially randomly distributed over the search space ($\Omega$). The particles move towards a global minima in each iteration and the direction of movement is given by the vector addition of the best value found so far (global best or *gbest*) among all the particles and the individual particle best position (*pbest*). The objective function which is essentially the minimization criteria is used to evaluate the fitness of the particle at a specific position in each iteration. For each particle ($i$) the velocity in each dimension in the next iteration is updated by the following velocity and position update equation, given by

$$
\begin{aligned}
v_i(t+1) &= \omega v_i(t) + c_1\phi_1(p_i(t) - x_i(t)) + c_2\phi_2(p_g(t) - x_i(t)) \\
x_i(t+1) &= x_i(t) + v_i(t+1)
\end{aligned}
\quad (16)
$$

Each particle's position ($x_i$) in the consequent iteration depends on its velocity ($v_i$) in the present iteration multiplied by an inertia factor ($\omega$) which is generally chosen large so as to prevent random erratic movement of the particles in the search space and to deviate the velocity of the particles by a smaller amount in each iteration. The other two positive constants $c_1, c_2$ denote the cognitive learning rate and the social learning rate respectively. The weights $c_1, c_2$ indicate the relative importance of the learning of the particles from its own best position ($p_i$) and the global best position ($p_g$) and both have been chosen as 1.4962 for the present study (Evers and Ben Ghalia 2009). In (16), $\{\phi_1, \phi_2\} \in [0,1]$ are two uniformly distributed random numbers. The inertia factor ($\omega$) is chosen to be 0.71633. A velocity clamping method is also incorporated in the algorithm and the maximum value of the velocity is set to 15% of the range in each dimension. This ensures that the velocity does not explode to bigger values and helps in controlling the global exploration of the swarm. The range here refers to the difference between the upper and the lower bounds of the search space ($\Omega$) on each dimension.

$$v_j^{max} = \lambda \cdot \text{range}_j(\Omega)$$
$$\text{range}_j(\Omega) = x_j^U - x_j^L, \ j = 1, 2, \cdots, n. \tag{17}$$

Here, $\text{range}_j(\Omega)$ represents range of search space along dimension $j$ and $\lambda$ represents the velocity clamping percentage. In the present case, the number of particles $N_p$ is chosen to be 20.

However the classical *gbest* PSO might get trapped in local minima resulting in premature convergence. This gives rise to stagnation where the particles continue to converge within a small region resulting in the fact that the global best and the personal bests are all within a very small radius in the search space. Thus the particles stagnate as the momentum from their previous velocities die out.

Regrouping PSO is a stochastic optimization algorithm which is an improvisation over the classical *gbest* PSO to prevent the problem of premature convergence and stagnation (Evers and Ben Ghalia 2009). This is done by an automatic triggering of the regrouping mechanism when premature convergence is detected. This method liberates the particles from these local minima and ensures convergence of the solution towards the true global minima. The swarm radius $\delta(t)$ at each iteration $t$ is taken as an indicator for premature convergence and gives the maximum Euclidean distance in the n-dimensional search space of any particle from the global best ($g(t)$). It is given by the following expression:

$$\delta(t) = \max_{\forall i} \|x_i(t) - g(t)\| \tag{18}$$

where, $\|\cdot\|$ denotes the Euclidean norm.

The regrouping mechanism is triggered when the normalized swarm radius ($\delta_{norm}$) is less than a pre-specified stagnation threshold ($\varepsilon$)

$$\text{i.e. } \delta_{norm} = \frac{\delta(t)}{\|range(\Omega)\|} \tag{19}$$

As in (Evers and Ben Ghalia 2009) ε is chosen as $1.1 \times 10^{-4}$ which is seen to work well with the regrouping mechanism for the optimization of wide range of problems. On detection of stagnation the swarm is regrouped in an area centered about the global best so that it is small enough for efficient search and at the same time ensures that the solution can go out of the local minima. The regrouping factor ($\rho$) is given by

$$\rho = \frac{6}{5\varepsilon} \tag{20}$$

which works well for all benchmark problems as indicated in (Evers and Ben Ghalia 2009).

When premature convergence is detected, the range in which particles are to be regrouped, around the global best is calculated for each dimension as the minimum of the original range of the search space on that dimension and the product of the regrouping factor with the maximum distance along that dimension of any particle from global best. i.e.

$$range_j(\Omega^r) = \max\left(range_j(\Omega^0), \rho \max_{i \in \{1,\cdots,s\}} |x_{i,j}^{r-1} - g_j^{r-1}|\right) \tag{21}$$

After re-initializing the particles' positions the swarm is now regrouped as

$$\begin{aligned} x_i &= p_g^{r-1} + r' \circ range(\Omega^r) - \frac{1}{2}range(\Omega^r) \\ range(\Omega^r) &= \left[range_1(\Omega^r), \cdots, range_n(\Omega^r)\right] \end{aligned} \tag{22}$$

utilizing a random vector $r'$ within the implicitly defined search space

$$\Omega^r = \left[x_1^{L,r}, x_1^{U,r}\right] \times \left[x_2^{L,r}, x_2^{U,r}\right] \times \cdots \times \left[x_n^{L,r}, x_n^{U,r}\right] \tag{23}$$

Having respective lower and upper bounds as

$$\begin{aligned} x_j^{L,r} &= p_{gj}^{r-1} - \frac{1}{2}range_j(\Omega^r) \\ x_j^{U,r} &= p_{gj}^{r-1} + \frac{1}{2}range_j(\Omega^r) \end{aligned} \tag{24}$$

Here, the operator " ∘ " denotes Hadamard element wise vector product, $r$ denotes swarm regrouping index, $p_g^{r-1}$ and $x_i^{r-1}$ denotes the global best and position of the i[th] particle at the last iteration of prior regrouping. For further details of regrouping PSO algorithm please refer to (Evers and Ben Ghalia 2009).

It is to be noted that the GA is utilized here in the inner loop and the PSO in the outer loop and not the other way round. It is well known that the PSO algorithm is a much better global optimizer than GA (Hassan et al. 2005; Ou and Lin 2006), for similar controller parameter searching purpose using time domain performance index (Pan et al. 2011b). Hence even if the GA is not able to find whether a particular solution is stable or not, the global search for the best solution still continues due to the regrouping nature of the PSO algorithm. A RegPSO-RegPSO algorithm could also have been used for both the stages, but it would have made the complexity of the problem unnecessarily higher while marginally affecting the outcome of the optimal solutions. This is because, the inner RegPSO would spend more time and computational resources on trying to find whether a particular solution is stable or not. Even if some solutions are not found to be stable, the outer loop would rely on the obtained stable solutions to direct the search towards finding optimal solutions which have good closed loop performance. The time spent refining solutions for just stability is not really as important as exploring other solutions affecting the performance of the NCS by the outer RegPSO loop. For a fair comparison, the convergence characteristics of the proposed RegPSO-GA hybrid method have been compared with the GA-GA based one as a part of the simulation example.

## 4. Illustrative example

A MATLAB/Simulink model has been developed to solve the problem and it has been run several times for each stable guess value of the state feedback controller within the optimization algorithm. Then the average of these is taken so that the expectation of global minima for the stochastic fitness function can be found out numerically.

A continuous time system (25) has been discretized with sampling time $h = 0.3$ for case study as reported in (Zhang et al. 2001). The packet drop-out probability has been considered as 30% for the present simulation study.

$$\begin{bmatrix} \dot{x}_1 \\ \dot{x}_2 \end{bmatrix} = \begin{bmatrix} 0 & 1 \\ 0 & -0.1 \end{bmatrix} \begin{bmatrix} x_1 \\ x_2 \end{bmatrix} + \begin{bmatrix} 0 \\ 0.1 \end{bmatrix} u,$$
$$y = \begin{bmatrix} 1 & 0 \end{bmatrix} \begin{bmatrix} x_1 \\ x_2 \end{bmatrix}. \tag{25}$$

For the discrete system (25) as reported in (Zhang et al. 2001), the proposed methodology is applied using the GA-GA hybrid optimization technique and the optimum values for the weighting matrices are reported below:

$$Q = \begin{bmatrix} 0.29495 & 0 \\ 0 & 1.37137 \end{bmatrix}, R = 0.25781 \text{ with minimized cost of } J_{min} = 111.4242.$$

These weighting matrices produce the sub-optimal stabilizing state-feedback gains as $K = \begin{bmatrix} 1.00337 & 4.09011 \end{bmatrix}$. The numerical values of the GA based BMI solver which gives the stability condition in the presence of packet losses are $a_1 = 1.0604, a_2 = 0.8772$. The $\tilde{P}$ matrix of the quadratic Lyapunov function has been computed as:

$$\tilde{P} = \begin{bmatrix} 27.779 & 39.8478 & -17.1649 & -14.5208 \\ 39.8478 & 188.0877 & -13.848 & -83.083 \\ -17.1649 & -13.848 & 17.5084 & 15.2433 \\ -14.5208 & -83.083 & 15.2433 & 87.2639 \end{bmatrix}$$

The eigen-values of the $\tilde{P}$ matrix is given as $eig(\tilde{P}) = \begin{bmatrix} 1.9366 \\ 30.6368 \\ 42.8462 \\ 245.2193 \end{bmatrix}$

Since the eigen-values of $\tilde{P}$ matrix are all greater than zero, the condition of positive semi-definiteness is satisfied. Equation (5) is also satisfied with these values of $a_1, a_2$ since the left hand side evaluates to $a_1^r a_2^{1-r} = 1.0017$ which is greater than unity.

The hybrid RegPSO-GA algorithm gives the optimized LQR weights as $Q = \begin{bmatrix} 11.87689 & 0 \\ 0 & 14.33702 \end{bmatrix}, R = 10.58286$, with a minimized expected cost of $J_{min} = 62.155$. It can be seen that the RegPSO-GA based algorithm gives a lower value of the expected objective function than the GA-GA based algorithm shown above. These weighting matrices produce the sub-optimal stabilizing state-feedback gains as $K = [0.99994 \quad 3.73058]$. The numerical values of the GA based BMI solver which gives the stability condition in the presence of packet losses are $a_1 = 1.0655, a_2 = 0.86331$. The $\tilde{P}$ matrix of the quadratic Lyapunov function has been computed as:

$$\tilde{P} = \begin{bmatrix} 11.7672 & 15.0175 & -7.7426 & -4.1701 \\ 15.0175 & 72.6810 & -5.1676 & -27.2520 \\ -7.7426 & -5.1676 & 7.6440 & 4.2129 \\ -4.1701 & -27.2520 & 4.2129 & 25.2804 \end{bmatrix}$$

The eigen-values of the $\tilde{P}$ matrix is given as $eig(\tilde{P}) = \begin{bmatrix} 0.7439 \\ 13.1253 \\ 14.5415 \\ 88.9619 \end{bmatrix}$

Since the eigen-values of $\tilde{P}$ matrix are all greater than zero, the condition of positive semi-definiteness is satisfied. Equation (5) is also satisfied with these values of $a_1, a_2$ since the left hand side evaluates to $a_1^r a_2^{1-r} = 1.0004$ which is greater than unity.

The convergence curves for the hybrid PSO-GA and hybrid GA-GA have been compared in Fig. 4. The convergence curves clearly show that the proposed method is capable of finding the global minima in much lesser number of generations. The comparative improvements are shown in Fig. 5, as the state trajectories and output of the system with packet losses with the GA-GA and RegPSO-GA hybrid algorithms which shows the later produces much faster set-point tracking performance even in the presence of random packet losses. In Fig. 4, the set-point tracking performances are compared for one single realization of the random variable representing the packet loss in the network, for the sake of simplicity, although the hybrid RegPSO-GA algorithm ensures minimization of the expected or average value of the set-point tracking error. Since, the performance comparison for such a stochastic process cannot be done from a single run, the proposed hybrid algorithms have been run multiple times and the convergence characteristics of corresponding to the lowest minima have been shown in Fig. 4.

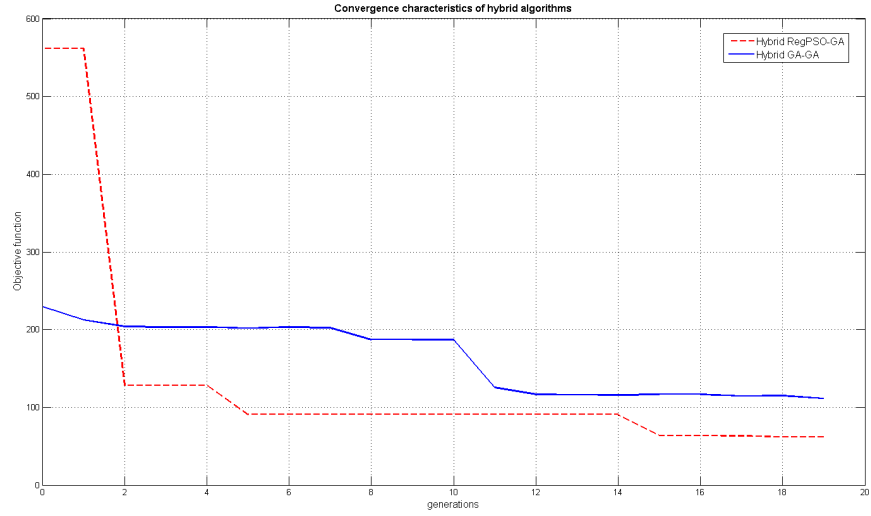

Fig. 4. Comparison of the convergence characteristics for hybrid RegPSO-GA and GA-GA algorithms.

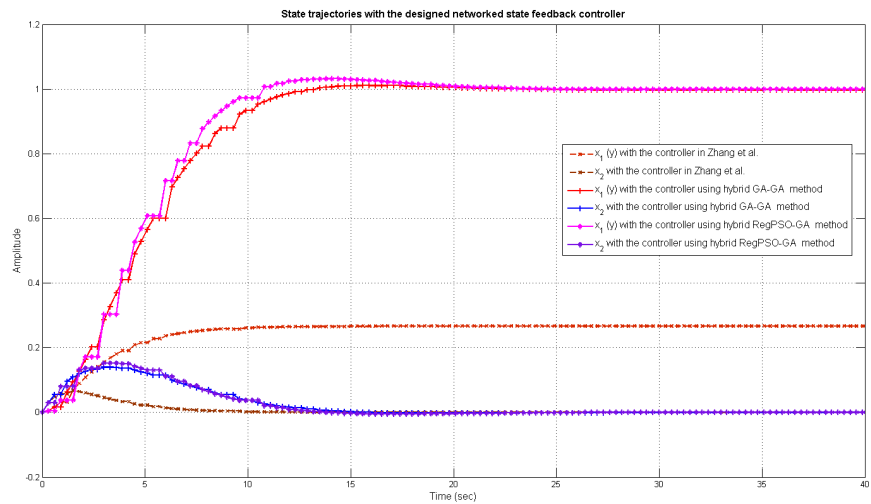

Fig. 5. State trajectories with the designed stabilizing sub-optimal controller with hybrid GA-GA and RegPSO-GA algorithm.

In (Zhang et al. 2001), the scaled output response is shown but in the actual output there is a steady state DC offset with the state feedback controller which is not desirable, especially in networked process control applications. Fig 5 shows that in the present method the steady state error is removed by the state feedback gains itself due to the inclusion of the additional performance criteria in the optimization algorithm. The GA has been used to find a feasible solution to the BMI problem. This is less susceptible to be trapped in local minima as finding only one feasible solution suffices. The RegPSO is a more robust solution searching scheme than the GA and has been used in the outer loop so that the optimal controller gains can be found which minimizes the objective function as stated in (15). It is to be noted that the output or the first state variable always tracks the input step excitation even in the presence of random packet loss as a stochastic phenomenon in the networked control loop. Also, since packet drop is a stochastic phenomenon in NCS, sudden jumps may be there in the state trajectories. For example see in (Pan et al. 2011a) where the state variables show sharp changes to compensate for the packet loss and to guarantee closed loop stability. The proposed technique ensures good set-point tracking performance which is evident from the simulation results.

The mechanism of simulating the packet drop has been incorporated using a MATLAB/Simulink based environment. The packet drop phenomenon in the feedback path as in Fig. 1 has a normal distribution. To illustrate further on the stochastic nature of the NCS, time domain simulation results for different realization of the random packet loss in the feedback path has been shown in Fig. 6 and Fig. 7 for the hybrid GA-GA and RegPSO-GA based state feedback controllers respectively. The corresponding values of the $J_{min}$ are also shown in the figures, which justifies that the proposed technique efficiently minimizes the expected value of the time domain performance index with guaranteed Lyapunov stability.

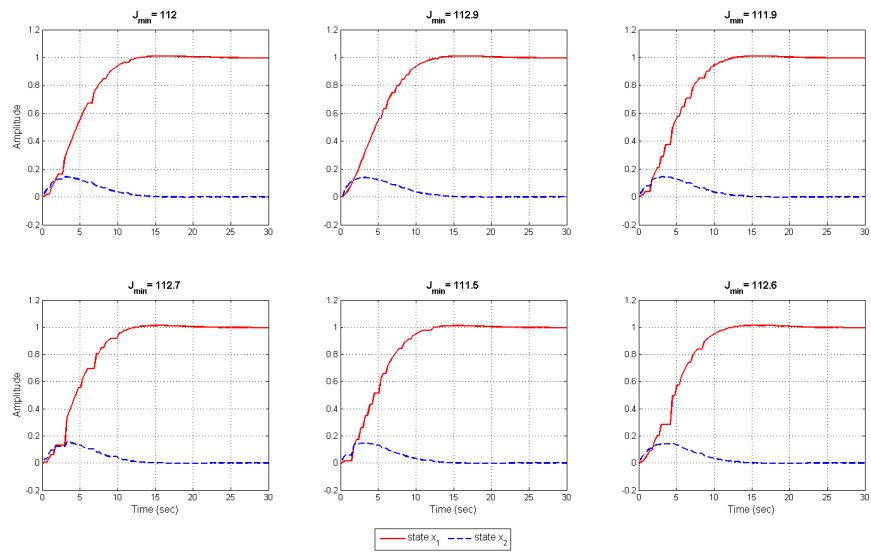

Fig. 6. State trajectories with the GA-GA hybrid algorithm based controller for different realization of the random packet loss phenomena.

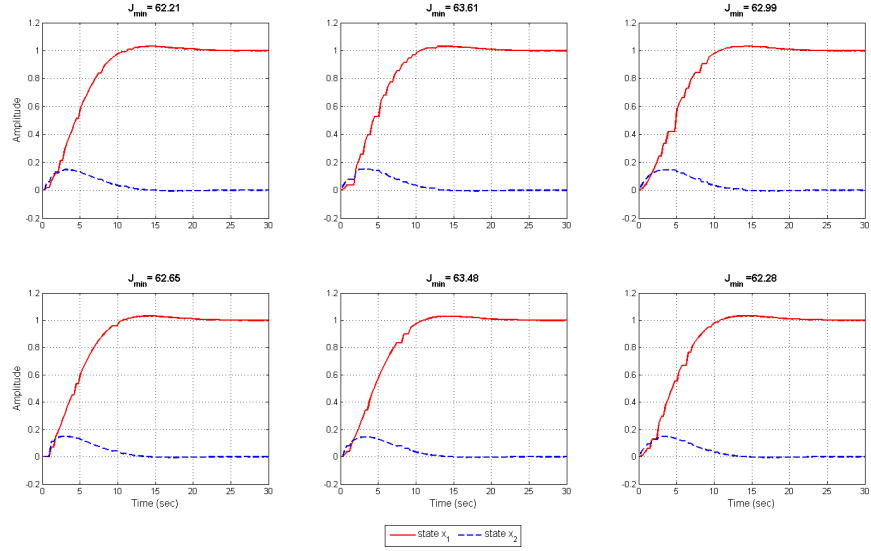

Fig. 7. State trajectories with the RegPSO-GA hybrid algorithm based controller for different realization of the random packet loss phenomena.

At a first glance, the time response curves for the PSO-GA (Fig. 7) and GA-GA (Fig. 6) cases seems to be closer but confusion may arise since the corresponding minima of the objective function is quite different. To illustrate this point more clearly a more detailed attention towards Fig. 5 is needed. It can be seen that the output state $x_1$ has slightly faster rise time for the RegPSO-GA algorithm than the GA-GA algorithm. Also, in equation (15) small variation in the state trajectories is continuously weighted by the time index ($k$). This basically translates to the fact that even small deviations from the best solution are given heavy penalty due to the time multiplication term (i.e. as time increases the penalty at the later stages for the system trajectory deviation is much higher). Thus even though Fig 5 shows a small improvement in performance visually, the time weighted penalty term assigns a lot of penalty to even slightly worse solutions. It is to be noted that the formulation could have been done without the time dependent term using an IAE (Integral of Absolute Error) criterion instead of the ITAE. This would have made the $J_{min}$ values nearly of the same magnitude with the PSO-GA being slightly better than the GA-GA. However this is not good from the control system design viewpoint as it is desirable for a system to reach steady state quickly and not have persistent oscillations at later time instants. Since, the IAE criteria would equally weigh the deviations from the set-point at all instants of time, the solutions which have a small oscillatory response at the later time instants can also be given lower fitness values if they have a fast rise time. From the perspective of the evolutionary/swarm computation, this ITAE criterion gives a higher selection pressure so that better solutions may be found.

The convergence curves shown in Fig. 4 correspond to the best found minima for the two algorithms among multiple runs. It is well known that comparison of statistical measures of two optimization algorithms are generally done to find out which algorithm has better consistency in finding the global minima. Thus the mean, standard-deviation, best and worst value of the expected objective function ($J_{\min}$) for 30 independent runs have been shown in Table 1 which shows that each of the statistical measures are better for RegPSO-GA algorithm compared to the GA-GA algorithm. It is to be noted that the motivation of applying GA/PSO to the current problem is different compared to the conventional optimization based controller design related literatures. The GA is applied to convert the BMI problem to be solvable using standard LMI solvers, though other complicated algorithms could have also been used which takes more computational resources. The PSO is applied to locate the global minima of the time domain tracking performance index. It is shown in (Pan et al. 2011b) that PSO has higher capability of finding the expected minima of tracking error in similar NCS design problems than that with GA. Apart from good set-point tracking with the consideration of the randomness of NCS, the algorithm also ensures Lyapunov stability by satisfying the BMI criterion using hybridization of standard LMI and GA solvers.

Table 1: Comparison of statistical performances of two hybrid algorithms for 30 independent runs.

| Hybrid Algorithm | Mean of $J_{\min}$ | Standard deviation of $J_{\min}$ | Best $J_{\min}$ | Worst $J_{\min}$ |
|---|---|---|---|---|
| GA-GA | 112.5867 | 0.7771 | 111.4 | 115 |
| PSO-GA | 63.1563 | 0.6573 | 62.15 | 64.54 |

The specific improvisations of the present work over state of the art literatures are summarized as follows.
- In (Zhang et al. 2001) the state feedback gains are pre-specified and the BMI problem is solved by a random search method. In the present paper optimization is done to get the state feedback controller gains to minimize the discrete time ITAE like criterion while simultaneously ensuring stability of the system by solution of the BMI problem. The GA based BMI solver which involves solution through both GA and interior point method simultaneously is an improvisation over (Sekaj and Vesely 2005) which becomes cumbersome when the number of problem dimensions of the GA increases due to increase in size of the solution matrix. The interior point methods can efficiently handle LMIs of large dimension and hence makes our solution methodology capable of handling solution variables of more dimensions.

- Unlike (Zhang et al. 2001), the choice of gains is obtained from the solution of the Discrete Algebraic Riccati Equation and hence deviation in the state variable trajectories and the control signal is minimized in an optimal fashion for the case without packet drops. In the presence of packet drops the solution becomes slightly suboptimal as the LQR controller is not optimal in the presence of packet losses in the network.

- In (Zhang et al. 2001), steady state offset is not eliminated, but in our case, the additional performance criteria (15) ensures that the steady state error is minimized along with a faster rise time, faster settling time and a lower peak overshoot.

- The weighting matrices of discrete time LQR (i.e. $Q$ and $R$) are optimally chosen with the RegPSO algorithm, as these are difficult to pre-specify at the beginning of the solution to ensure set-point tracking even in the presence of packet losses. In classical LQR, the weights are chosen a priori and also LQR formulation can't guarantee closed loop stability in the presence of packet losses in the network. Here optimum choice of weighting matrices of LQR is done while also ensuring good set-point tracking performance. Guaranteed Lyapunov stability comes from the BMI formulation itself for arbitrary packet drops in the network. Traditional LQR neither can guarantee stability for random packet losses nor is capable of tracking step reference input for arbitrary weighting matrices. Herein lies the motivation of the present paper, with RegPSO based design of optimum LQR weights while solving the BMI problem by GA and standard LMI solvers.

- In this method only discrete time ITAE like performance index has been used as an additional control system performance criterion, but the fitness function can also be designed to handle other performance criteria like higher moments of time and tracking error along with the total variation of the control signal etc. (Pan et al. 2011b) with user specified weights attached to each criterion.

- In (Pan et al. 2011b) GA/PSO based optimization framework has been proposed to design fuzzy PID controllers for NCS applications. Though such an optimization minimizes the expected cost function involving the tracking performance of the system's output, the analytical stability criterion for fuzzy controllers over NCS is difficult to derive. In the present paper, with a relatively simpler state feedback controller, the NCS is modeled as a switched system and the analytical Lyapunov stability conditions are ensured which leads to the BMI problem. The BMI problem is then solved using the hybridization of LMI and GA, while the RegPSO minimizes the expectation of the performance index relating the tracking error. Thus the current paper improvises the concept of stochastic optimization based networked controller design with guaranteed analytical Lyapunov stability condition incorporated as a BMI problem within the optimization for closed loop performance.

It is an obvious fact that such hybrid algorithms are often questioned for the probable increase in computational complexity. Here, the computational complexity increases due to hybrid PSO-GA algorithms, along with the LQR with optimum weight selection. The achievable goal of (Zhang et al. 2001) was only stabilization of NCS and those of the present paper are optimum weight based LQR for set-point tracking in presence of packet losses in NCS. It is true that so many complex algorithms for performance improvement of a robust tracking problem in NCS would increase the complexity. But since the computation of the controller gains is offline, hence the additional complexity to increase the performance does not pose a significant problem.

## 5. Conclusion

A practical solution for controller designing in NCS applications in the presence of packet losses have been proposed in this paper using a hybrid Regrouping PSO-GA based algorithm. The stability of the closed loop system is guaranteed by Lyapunov based BMIs which is incorporated in the hybrid optimization algorithm itself. The consequent optimally chosen weighting matrices based discrete time LQR controller designing methodology with additional design constraints as the guaranteed stabilizing sub-optimal state feedback gain selection are shown by means of an illustrative example. Further work can be directed towards similar NCS design involving multi-loop systems.

## Acknowledgement

The authors thank the anonymous reviewers for providing helpful and constructive comments which has helped to increase the quality of the paper.

## References

Doyle JC (1983) Synthesis of robust controllers and filters. Conference on Decision and Control, The 2nd IEEE Conference on, San Antonio, TX 109–114.

El Ghaoui L, Balakrishnan V (1994) Synthesis of fixed-structure controllers via numerical optimization. Decision and Control, 1994, Proceedings of the 33rd IEEE Conference on 3:2678–2683.

Evers GI, Ben Ghalia M (2009) Regrouping Particle Swarm Optimization: A New Global Optimization Algorithm with Improved Performance Consistency Across Benchmarks. Systems, Man and Cybernetics, 2009 SMC 2009 IEEE International Conference on 3901–3908.

Fukuda M, Kojima M (2001) Branch-and-cut algorithms for the bilinear matrix inequality eigenvalue problem. Computational Optimization and Applications 19:79–105.

Goh KC, Safonov MG, Papavassilopoulos GP (1995) Global optimization for the biaffine matrix inequality problem. Journal of global optimization 7:365–380.

Hassan R, Cohanim B, De Weck O, Venter G (2005) A comparison of particle swarm optimization and the genetic algorithm. Proceedings of the 1st AIAA Multidisciplinary Design Optimization Specialist Conference

Iwasaki T (1999) The dual iteration for fixed-order control. Automatic Control, IEEE Transactions on 44:783–788.

Kanev S, Scherer C, Verhaegen M, De Schutter B (2004) Robust output-feedback controller design via local BMI optimization. Automatica 40:1115–1127.

Kennedy J, Eberhart R (1995) Particle swarm optimization. Neural Networks, 1995 Proceedings, IEEE International Conference on 4:1942–1948.

Kim E, Kim S (2002) Stability analysis and synthesis for an affine fuzzy control system via LMI and ILMI: a continuous case. Fuzzy Systems, IEEE Transactions on 10:391–400.

Kong Y, Zhao D, Yang B, et al. (2012) Static output feedback control for active suspension using PSO-DE/LMI approach. Mechatronics and Automation (ICMA), 2012 International Conference on 366–370.

Lo JC, Lin ML (2004) Observer-based robust H∞ control for fuzzy systems using two-step procedure. Fuzzy Systems, IEEE Transactions on 12:350–359.

Longo S, Herrmann G, Barber P (2012) Robust Scheduling of Sampled-Data Networked Control Systems. Control Systems Technology, IEEE Transactions on 20:1613–1621.

Longo S, Su T, Herrmann G, et al. (2011) Scheduling of the FlexRay static segment for robust controller integration. Control Applications (CCA), 2011 IEEE International Conference on 1487–1492.

Ogata K (1995) Discrete-time control systems. Prentice-Hall Englewood Cliffs, NJ

Ou C, Lin W (2006) Comparison between PSO and GA for Parameters Optimization of PID Controller. Mechatronics and Automation, Proceedings of the 2006 IEEE International Conference on 2471–2475.

Pan I, Das S, Ghosh S, Gupta A (2011a) Stabilizing Gain Selection of Networked Variable Gain Controller to Maximize Robustness Using Particle Swarm Optimization. Process Automation, Control and Computing (PACC), 2011 International Conference on 1–6.

Pan I, Das S, Gupta A (2011b) Tuning of an optimal fuzzy PID controller with stochastic algorithms for networked control systems with random time delay. ISA Transactions 50:28–36.

Sekaj I, Vesely V (2005) Robust output feedback controller design: genetic algorithm approach. IMA Journal of Mathematical Control and Information 22:257–265.

Tipsuwan Y, Chow MY (2003) Control methodologies in networked control systems. Control engineering practice 11:1099–1111.

Toker O, Ozbay H (1995) On the NP-hardness of solving bilinear matrix inequalities and simultaneous stabilization with static output feedback. American Control Conference, 1995 Proceedings of the 4:2525–2526.

Tuan H, Apkarian P (2000) Low nonconvexity-rank bilinear matrix inequalities: algorithms and applications in robust controller and structure designs. Automatic Control, IEEE Transactions on 45:2111–2117.

Tuan H, Apkarian P, Narikiyo T, Yamamoto Y (2001) Parameterized linear matrix inequality techniques in fuzzy control system design. Fuzzy Systems, IEEE Transactions on 9:324–332.

VanAntwerp JG, Braatz RD (2000) A tutorial on linear and bilinear matrix inequalities. Journal of Process Control 10:363–385.

Vandenberghe L, Balakrishnan V (1997) Algorithms and software for LMI problems in control. Control Systems, IEEE 17:89–95.

Yamada Y, Hara S (1998) Global optimization for H∞ control with constant diagonal scaling. Automatic Control, IEEE Transactions on 43:191–203.

Zamani M, Sadati N, Ghartemani MK (2009) Design of an H∞ PID controller using particle swarm optimization. International Journal of Control, Automation and Systems 7:273–280.

Zhang W, Branicky MS, Phillips SM (2001) Stability of networked control systems. Control Systems, IEEE 21:84–99.